\documentclass[leqno,11pt]{article}  
\usepackage{amsmath}
\usepackage{amssymb}

\usepackage[german,english]{babel}
\usepackage{amsthm}

\title{Acyclic coefficient systems on buildings}

\author{\textsc{Elmar Grosse-Kl\"onne}}
\date{}

\theoremstyle{plain} 
\newtheorem{satz}{Theorem}[section]  
\newtheorem{kor}[satz]{Corollary}  
\newtheorem{lem}[satz]{Lemma}  
\newtheorem{pro}[satz]{Proposition}  
 
\newcommand{\bi}{\mbox{\rm im}}  
\newcommand{\ke}{\mbox{\rm Ker}}  

\newcommand{\kara}{\mbox{\rm char}}  

\newcommand{\id}{\mbox{\rm id}}

\newcommand{\rk}{\mbox{\rm rank}}

\newcommand{\Gr}{\mbox{\rm Gr}}

\theoremstyle{remark}

\theoremstyle{definition}

\begin{document}
\maketitle

\footnote[0]{I am very grateful to Peter Schneider whose comments helped me a lot to eliminate inacurracies and to streamline presentation, terminology and notations. I thank Ehud de Shalit and Annette Werner for related discussions. Thanks also go to the referee for his careful reading and for supplying a simple proof for Lemma \ref{collex}.}

\begin{abstract}  For cohomological (resp. homological) coefficient systems ${\mathcal F}$ (resp. ${\mathcal V}$) on affine buildings $X$ with Coxeter data of type $\widetilde{A}_d$ we give for any $k\ge1$ a sufficient local criterion which implies $H^k(X,{\mathcal F})=0$ (resp. $H_k(X,{\mathcal V})=0)$. Using this criterion we prove a conjecture of de Shalit on the acyclicity of coefficient systems attached to hyperplane arrangements on the Bruhat-Tits building of the general linear group over a local field. We also generalize an acyclicity theorem of Schneider and Stuhler on coefficient systems attached to representations. 
 
\end{abstract}

%


\begin{center} {\bf Introduction} \end{center}

Let $X$ be an affine building whose apartments are Coxeter complexes attached to Coxeter systems of type $\widetilde{A}_d$ and let ${\mathcal F}$ be a cohomological coefficient system (ccs) on $X$. The purpose of this paper is to give a {\it local} criterion which assures that for a given $k\ge1$ the cohomology group $H^k(X,\mathcal F)$ vanishes. Similarly for homological coefficient systems.

For a sheaf ${\mathcal G}$ on a topological space $Y$ it is well known that $H^k(Y,\mathcal G)=0$ for all $k\ge1$ if ${\mathcal G}$ is {\it flasque}, i.e. if all restriction maps ${\mathcal G}(U)\to {\mathcal G}(V)$ for open $V\subset U\subset Y$ are surjective. We are looking for an adequate notion of "flasque"\, ccs on $X$.

If $d=1$ the same condition works: If the restriction map ${\mathcal F}(\sigma)\to{\mathcal F}(\tau)$ for any $0$-simplex $\sigma$ contained in the $1$-simplex $\tau$ is surjective, then $H^1(X,{\mathcal F})=0$. This is easily seen using the {\it contractibility} of $X$. However, if $d>1$ the surjectivity of ${\mathcal F}(\sigma)\to{\mathcal F}(\tau)$ for any $(k-1)$-simplex $\sigma$ contained in the $k$-simplex $\tau$ does not guarantee $H^k(X,{\mathcal F})=0$. The other naive transposition of the flasqueness concept from topological spaces to buildings would be to require for any $(k-1)$-simplex $\sigma$ the surjectivity of ${\mathcal F}(\sigma)\to\prod_{\tau}{\mathcal F}(\tau)$, taking the product over {\it all} $k$-simplices $\tau$ containing $\sigma$. This would indeed force $H^k(X,{\mathcal F})=0$ but would also be a completely useless criterion: for example, it would not be satisfied by a constant ccs ${\mathcal F}$ (of which we know $H^k(X,{\mathcal F})=0$, by the contractibility of $X$).

Let us describe our criterion ${\mathcal S}(k)$. We fix an orientation of $X$. It defines a {\it cyclic} ordering on the set of vertices of any simplex, hence a true ordering on the set of vertices of any {\it pointed} simplex. To a pointed $(k-1)$-simplex $\widehat{\eta}$ we associate the set $N_{\widehat{\eta}}$ of all vertices $z$ for which $(z,\widehat{\eta})$ (i.e. $z$ as the first vertex) is an ordered $k$-simplex (in the previously qualified sense). We define what it means for a subset $M_0$ of $N_{\widehat{\eta}}$ to be {\it stable with respect to $\widehat{\eta}$} (if for example $d=1$ the condition is $|M_0|\le1$). Our criterion ${\mathcal S}(k)$ which assures $H^k(X,{\mathcal F})=0$ is then that for any such $\widehat{\eta}$ and for any subset $M_0$ of $N_{\widehat{\eta}}$, stable with respect to $\widehat{\eta}$, the sequence$${\mathcal F}(\eta)\longrightarrow\prod_{z\in {M}_0}{\mathcal 
F}(\{z\}\cup\eta)\longrightarrow\prod_{z,z'\in {M_0}}{\mathcal 
F}(\{z,z'\}\cup\eta)$$is exact (where in the target of the second arrow the product is over pairs $z,z'\in M_0$ of incident vertices). For example, a constant ccs ${\mathcal F}$ satisfies ${\mathcal S}(k)$. 

Having fixed a vertex $z_0$, the central ingredient in the proof is then a certain function $i$ on the set $X^0$ of vertices which measures the combinatorical (not: Euclidean) distance from $z_0$; it depends on our chosen orientation.

Dual to ${\mathcal S}(k)$ we describe a  criterion ${\mathcal S}^*(k)$ which guarantees $H_k(X,{\mathcal V})=0$ for homological coefficient systems (hcs) ${\mathcal V}$ on $X$.

The basic example is of course the case where $X$ is the Bruhat-Tits building of ${\rm PGL_{d+1}}(K)$ for a local field $K$. Our original motivation for developing the criterion ${\mathcal S}(k)$ was the following. In \cite{schn}, Schneider defined a certain class of ${\rm SL_{d+1}}(K)$-representations, the "holomorphic discrete series representations", as the global sections of certain equivariant vector bundles on  Drinfel'd's symmetric space $\Omega_{K}^{(d+1)}$ of dimension $d$ over $K$. In a subsequent paper \cite{int} we will construct for any such vector bundle $V$ an integral model ${\mathcal V}$, as an equivariant coherent sheaf on the formal ${\mathcal O}_K$-scheme $\Omega_{{\mathcal O}_K}^{(d+1)}$ underlying $\Omega_{K}^{(d+1)}$. Using the criterion ${\mathcal S}(k)$ and the close relation between $\Omega_{K}^{(d+1)}$ and $X$, we will show that if $V$ is strongly dominant  (in a suitable sense), then $H^k(\Omega_{{\mathcal O}_K}^{(d+1)},{\mathcal V})=0$ for all $k\ge1$. Examples for such ${\mathcal V}$ are the terms of the logarithmic de Rham complex of $\Omega_{{\mathcal O}_K}^{(d+1)}$.

Here we present two other applications. The first one is a proof of a conjecture of de Shalit on $p$-adic ccs of Orlik-Solomon algebras. For an arbitrary field $K$, the assignment of the Orlik-Solomon algebra $A$ to a hyperplane arrangement $W$ in $(K^{d+1})^*$ --- the complement of a {\it finite} set ${\mathcal A}\subset {\mathbb{P}}(K^{d+1})$ of hyperplanes in $(K^{d+1})^*$ --- is a classical theme. $A$ is defined combinatorically in terms of the hyperplanes and turns out to be isomorphic with the cohomology ring of $W$. Of course $K$ may also be a local field. However, de Shalit \cite{ds} discovered that one can go further and give the story a {\it genuinely $p$-adic} flavour. Namely, if $K$ is a local field he allows ${\mathcal A}\subset{\mathbb{P}}(K^{d+1})$ to be {\it infinite}. He assigns to ${\mathcal A}$ not a single Orlik-Solomon algebra but a ccs $A=A(.)$ of Orlik-Solomon algebras on the Bruhat-Tits building $X$ of ${\rm PGL_{d+1}}(K)$ which then should play the role of a cohomology ring of the "hyperplane arrangement"\, defined by ${\mathcal A}$. The algebra $A(\sigma)$ for a $j$-simplex $\sigma$ is closely related to a suitable tensor product of Orlik-Solomon algebras for finite hyperplane arrangements in $k$-vector spaces, for $k$ the residue field of $K$. de Shalit conjectures that these beautiful ccs are acyclic in positive degrees, $H^k(X,A)=0$ for all $k\ge1$. He proved the conjecture for any ${\mathcal A}$ if $d\le 2$. For arbitrary $d$ he proved it if ${\mathcal A}$ is the full set ${\mathbb{P}}(K^{d+1})$ of all $K$-rational hyperplanes, and Alon \cite{alon} proved it if ${\mathcal A}$ is finite. Here we give a proof for all ${\mathcal A}$ and $d$ by showing that $A$ always satisfies ${\mathcal S}(k)$. 

In fact we prove a version with arbitrary coefficient ring: While in \cite{ds} the coefficient ring of the ccs $A$ is $K$, we allow an arbitrary coefficient ring $R$, e.g. also $R=\mathbb{Z}$ or $R=k$. While de Shalit's proof in case $d\le 2$ also works for arbitrary $R$, his proof in case ${\mathcal A}={\mathbb{P}}(K^{d+1})$ but $d$ arbitrary, which is by reduction to the main result of \cite{ss}, does not work for coefficient rings $R$ other than characteristic zero fields. We explain why this improvement for ${\mathcal A}={\mathbb{P}}(K^{d+1})$ and $R=\mathbb{Z}$ should have an application to a problem on $p$-adic Abel-Jacobi mappings raised in \cite{raxa}. 

The second application we describe is concerned with the technique of Schneider and Stuhler to spread out representations of ${\rm GL_{d+1}}(K)$ as hcs on the Bruhat-Tits building of ${\rm PGL_{d+1}}(K)$. For a ${\rm GL_{d+1}}(K)$-representation on a (not necessarily free) $\mathbb{Z}[\frac{1}{p}]$-module $V$ (where $p=\kara(k)$) which is generated by its vectors fixed under a principal congruence subgroup $U^{(n)}$ of some level $n>1$, we prove that the chain complex of the corresponding hcs is a resolution of $V$. For fields of characteristic zero as coefficient ring (instead of $\mathbb{Z}[\frac{1}{p}]$) this is the main result of \cite{ss} --- where however $n=1$ is allowed. While the proof in \cite{ss} uses the Bernstein-Borel-Matsumoto theory, we do not need any representation theoretic input whatsoever.

\section{The criterion}

\label{thec}

Let $d\ge1$ and let $X$ be an affine building whose apartments are Coxeter complexes attached to Coxeter systems of type $\widetilde{A}_d$. We refer to the book \cite{bro} for the basic definitions and properties of buildings. For $0\le j\le d$ we denote by $X^j$ the set of $j$-simplices. We generally identify a $j$-simplex with its set of vertices.  

We fix an orientation of $X$. It distinguishes for any simplex a {\it cyclic} ordering on its set of vertices. A {\it pointed $k$-simplex} is an {\it enumeration} of the set of vertices of a $k$-simplex in its distinguished cyclic ordering; we write it as an ordered $(k+1)$-tuple of vertices. 

For an apartment $A$ in $X$ we will slightly abuse notation by not distinguishing between $A$ and its geometric realization $|A|$. There is (see \cite{bro} p.148) an isomorphism of $A$ with $\mathbb{R}^{d+1}/{\mathbb{R}}.(1,\ldots,1)$ --- we view it here as an identification --- such that, if $\{e_0,\ldots,e_d\}$ denotes the standard basis of $\mathbb{R}^{d+1}$ the following holds:\begin{itemize}\item the set of vertices in $A$ is $\mathbb{Z}^{d+1}/{\mathbb{Z}}.(1,\ldots,1)$\item a $k+1$-tuple $(x_0,\ldots,x_k)$ of vertices in $A$ is a pointed $k$-simplex if and only if there is a sequence $$\emptyset\ne J_0\subsetneq \ldots\subsetneq J_{k-1}\subsetneq\{0,\ldots,d\}$$ such that $\sum_{j\in J_t}{e}_j$ represents $x_t-x_k$ (formed with respect to the obvious group structure on $\mathbb{Z}^{d+1}/{\mathbb{Z}}.(1,\ldots,1)$), for any $0\le t \le k-1$.\end{itemize}

If $(x_0,\ldots,x_k)$ is a pointed $k$-simplex, we define $\ell((x_0,\ldots,x_k))$ to be the maximal number $r$ such that there exists a pointed $r$-simplex $(y_0,\ldots,y_r)$ with $x_0=y_0$ and $x_k=y_r$. For a pointed $(k-1)$-simplex $\widehat{\eta}=(x_1,\ldots,x_k)$ we define the set$$N_{\widehat{\eta}}=\{z\in X^0|\,\,(z,x_1,\ldots,x_k)\,\mbox{is a pointed}\,\,k\mbox{-simplex}\}.$$For $z\in N_{\widehat{\eta}}$ we write $(z,\widehat{\eta})$ for the pointed $k$-simplex $(z,x_1,\ldots,x_k)$. We define a partial ordering $\le$ on $N_{\widehat{\eta}}$ by$$u_1\le u_2\Longleftrightarrow[u_1=u_2\,\,\mbox{or}\,\,(u_1,u_2,\widehat{\eta})\,\,\mbox{is a pointed}\,\,(k+1)\mbox{-simplex}].$$

\begin{lem}\label{grund} For any $u_1, u_2\in N_{\widehat{\eta}}$ the set $$W^{\widehat{\eta}}_{u_1,u_2}=\{u\in N_{\widehat{\eta}}|\,\,u\le u_1\,\,\mbox{and}\,\,u\le u_2\}$$is empty or it contains an element $u$ such that $\ell((u,\widehat{\eta}))<\ell((w,\widehat{\eta}))$ for all $w\in W^{\widehat{\eta}}_{u_1,u_2}-\{u\}$.
\end{lem}

{\sc Proof:} Suppose we have two such candidates $u, u'$. We can find an apartment $A$ which contains $\eta$, $u_1$, $u_2$, $u$ and $u'$ (for example because $\eta\cup\{u_1,u\}$ and $\{u',u_2\}$ are simplices). We identify $A$ with $\mathbb{R}^{d+1}/{\mathbb{R}}.(1,\ldots,1)$ as above. There exist subsets $\emptyset\ne J_t\subsetneq\{0,\ldots,d\}$ for $1\le t\le k-1$ and $\emptyset\ne I_s\subsetneq\{0,\ldots,d\}$ for $s=1,2$ such that $x_t-x_k$ is represented by $\sum_{j\in J_t}e_j$ and $u_s-x_k$ is represented by $\sum_{j\in I_s}e_j$. We have $I_s\subsetneq J_1\subsetneq\ldots\subsetneq J_{k-1}$ for $s=1,2$. Hence both $u-x_k$ and $u'-x_k$ are represented by $\sum_{j\in I_1\cap I_2}e_j$. \\

If $W^{\widehat{\eta}}_{u_1,u_2}\ne \emptyset$ we denote the element $u\in W^{\widehat{\eta}}_{u_1,u_2}$ from \ref{grund} by $[\eta|u_1,u_2]$. If $W^{\widehat{\eta}}_{u_1,u_2}=\emptyset$ then $[\eta|u_1,u_2]$ is undefined. A subset $M_0$ of $N_{\widehat{\eta}}$ is called {\it stable with respect to $\widehat{\eta}$} if for any two vertices $u_1, u_2\in M_0$ the vertex $[\widehat{\eta}|u_1,u_2]$ is defined and belongs to $M_0$. (See \ref{fulat} below for what this means on Bruhat-Tits buildings. When working out the applications described later, we became fully convinced that stability is a very natural condition.) Note that $M_0$ is stable with respect to $\widehat{\eta}=(x_1,\ldots,x_k)$ if and only if it is stable with respect to the pointed $0$-simplex $x_k$ --- this is because of $[x_k|u_1,u_2]=[\eta|u_1,u_2]$, as we saw in the proof of \ref{grund}.\\
 
A {\it cohomological coefficient system} (ccs) ${\mathcal F}$ on $X$ is the assignment of an abelian group ${\mathcal F}({\tau})$ to every simplex ${\tau}$ of $X$, and a homomorphism $r^{{\tau}}_{{\sigma}}:{\mathcal F}({\tau})\to{\mathcal F}({\sigma})$ to every face inclusion ${\tau}\subset{\sigma}$, such that $r_{\rho}^{{\sigma}}\circ r_{{\sigma}}^{{\tau}}=r_{{\rho}}^{{\tau}}$ whenever ${\tau}\subset{\sigma}\subset{\rho}$, and $r_{{\tau}}^{{\tau}}$ is the identity.

Given a ccs ${\mathcal F}$, define the group $C^k(X,{\mathcal F})$ of $k$-cochains ($0\le k\le d$) to consist of the maps $c$, assigning to each $k$-simplex ${\tau}$ an element $c_{{\tau}}\in{\mathcal F}({\tau})$. Define$$\partial=\partial^{k+1}:C^k(X,{\mathcal F})\longrightarrow C^{k+1}(X,{\mathcal F})$$by the rule$$(\partial c)_{{\tau}}=\sum_{\tau'\subset\tau}[\tau:\tau']r_{{\tau}}^{{\tau}'}(c_{{\tau}'})$$where $[\tau:\tau']=\pm 1$ is the incidence number (with respect to a fixed labelling of $X$ as in \cite{bro} p. 30). Then $(C^{\bullet}(X,{\mathcal F}),\partial)$ is a complex ($\partial^2=0$), and its cohomology groups are denoted $H^k(X,{\mathcal F})$.\\
 
Consider for $1\le k\le d$ the following condition ${\mathcal S}(k)$ for a ccs ${\mathcal F}$ on $X$: For any pointed $(k-1)$-simplex $\widehat{\eta}$ with underlying $(k-1)$-simplex $\eta$ and for any subset $M_0$ of $N_{\widehat{\eta}}$ which is stable with respect to $\widehat{\eta}$, the following subquotient complex of $C^{\bullet}(X,{\mathcal F})$ is exact:$${\mathcal F}(\eta)\stackrel{\partial^k}{\longrightarrow}\prod_{z\in {M}_0}{\mathcal 
F}(\{z\}\cup\eta)\stackrel{\partial^{k+1}}{\longrightarrow}\prod_{z,z'\in {M_0}\atop\{z,z'\}\in{X^1}}{\mathcal 
F}(\{z,z'\}\cup\eta)$$(We regard the first term as a subgroup of $C^{k-1}(X,{\mathcal F})$, the second one as a direct summand of $C^{k}(X,{\mathcal F})$, and the third one as a quotient of $C^{k+1}(X,{\mathcal F})$.) Note that ${\mathcal S}(k)$ depends on the chosen orientation of $X$.

\begin{satz}\label{acykrit} Let ${\mathcal F}$ be a ccs on $X$. Let $1\le k\le d$ and suppose ${\mathcal S}(k)$ holds true. Then $H^k(X,{\mathcal F})=0$.  
\end{satz}

We fix once and for all a vertex $z_0\in X^0$. Given an arbitrary vertex $x\in X^0$, choose an apartment $A$ containing $z_0$ and $x$. Choose an identification of $A$ with $\mathbb{R}^{d+1}/{\mathbb{R}}.(1,\ldots,1)$ as before, but now require in addition that $z_0\in A$ corresponds to the class of the origin in $\mathbb{R}^{d+1}/{\mathbb{R}}.(1,\ldots,1)$. Let $\sum_{j=0}^dm_je_j$ be the unique representative of $x$ for which $m_j\ge0$ for all $j$, and $m_j=0$ for at least one $j$. Let $\pi$ be a permutation of $\{0,\ldots,d\}$ such that $m_{\pi(d)}\ge\ldots\ge{m}_{\pi(0)}\ge0$ and set$$i(x)=({m}_{\pi(d)},\ldots,{m}_{\pi(1)})\in\mathbb{N}_0^{d}.$$

\begin{lem} The $d$-tuple $i(x)$ is independent of the choice of $A$.
\end{lem}

{\sc Proof:} Let us write $i_A(x)$ instead of $i(x)$ in order to indicate the reference to $A$ in the above definition. Suppose the apartment $A'$ also contains $z_0$ and $x$. Choose a chamber ($d$-simplex) $C$ in $A$ containing $x$, and a chamber $C'$ in $A'$ containing $z_0$. Choose an apartment $A''$ in $X$ containing $C$ and $C'$, and let $\pi:X\to A''$, resp. $\pi':X\to A''$, be the retraction from $X$ to $A''$ centered in $C$, resp. centered in $C'$ (see \cite{bro} p.86). Then $\pi$, resp. $\pi'$, induces an isomorphism of oriented chamber complexes $A\cong A''$, resp. $A'\cong A''$. Hence $i_A(x)=i_{A''}(x)=i_{A'}(x)$.\\

Here is another, equivalent but more intrinsic definition of $i(x)$ (we do not need it). For $x,y\in X^0$ let $d(x,y)\in\mathbb{Z}_{\ge0}$ be the minimal number $t$ such that there exists a sequence $x_0,\ldots,x_t$ in $X^0$ with $x=x_0$, $y=x_t$ and $\{x_{r-1},x_{r}\}\in X^1$ for all $1\le r\le t$. For $x\in X^0$ and a subset $W\subset X^0$ let $d(x,W)=\min\{d(x,y)|\,y\in W\}$. For $1\le i\le d$ define the subset $W_i$ of $X^0$ inductively as follows: $W_1=\{z_0\}$ and $$W_i=\{z\in X^0|\,\,\left\{\begin{array}{c}\mbox{there exist elements}\,u_0,\ldots,u_r\in X^0\,\mbox{such that}\,r=d(z,W_{i-1}),\\u_0\in W_{i-1}, u_r=z\,\mbox{and}\,\ell((u_{\ell-1},u_{\ell}))=1\,\mbox{for all}\,1\le\ell\le r\end{array}\right\}\}.$$In particular $W_1\subset W_2\subset\ldots\subset W_d$. For $x\in X^0$ we then have$$i(x)=(d(x,W_1),d(x,W_2),\ldots,d(x,W_d)).$$Yet another equivalent definition of $i(x)$ (which we do not need either) results from the fact that the type of a minimal chamber-gallery connecting $x$ and $z_0$ encodes $i(x)$ if $x$ and $z_0$ are {\it not} incident.
 
On the set of ordered $d$-tuples $(n_1,\ldots,n_d)\in\mathbb{Z}_{\ge0}^d$ (and hence on the set of $d$-tuples $i(x)$ for $x\in X^0$) we use the lexicographical ordering: $$(n_1,\ldots,n_d)<(n'_1,\ldots,n'_d)\Longleftrightarrow\left\{\begin{array}{c}\mbox{there is a}\,\,1\le r\le d\,\,\mbox{such that}\\n_j=n'_j\,\,\mbox{for}\,\,1\le j\le r-1\,\,\mbox{and}\,\,n_r<n'_r\end{array}\right\}$$

\begin{lem}\label{iordpoin} Let ${\eta}$ be a $(k-1)$-simplex and let $x_1,\ldots,x_k$ be an enumeration of its vertices which satisfies $i(x_1)\le\ldots\le i(x_k)$. Then in fact $i(x_1)<\ldots<i(x_k)$ and $\widehat{\eta}=(x_1,\ldots,x_k)$ is a pointed $(k-1)$-simplex.
\end{lem}

{\sc Proof:} Choose an apartment $A$ containing $z_0$ and $\eta$, and choose an identification of $A$ with $\mathbb{R}^{d+1}/{\mathbb{R}}.(1,\ldots,1)$ as before, with $z_0\in A$ corresponding to the class of the origin in $\mathbb{R}^{d+1}/{\mathbb{R}}.(1,\ldots,1)$. Then the claims follow easily from our description of the simplicial structure of $\mathbb{Z}^{d+1}/{\mathbb{Z}}.(1,\ldots,1)$.\\ 

\begin{lem}\label{minnei} For any $x\in X^{0}$, $x\ne z_0$, there is among the vertices incident to $x$ a {\it unique} vertex $\nu(x)$ with minimal $i$-value: for all other vertices $z$ incident to $x$ we have $i(\nu(x))<i(z)$. If $z$ is incident to $x$, different from $\nu(x)$ and satisfies $i(z)<i(x)$, then $\nu(x)$ and $z$ are incident and $\ell((\nu(x),z))\le \ell((x,z))$. 
\end{lem}

{\sc Proof:} Let $[z_0,x]$ be the geodesic (with respect to the Euclidean distance function on the geometric realization $|X|$ of $X$) between $z_0$ and $x$. Let $\tau$ be the minimal simplex which contains $x$ and whose open interior (viewed as a subset of $|X|$) contains a point of $[z_0,x]$. We assert that the vertex $\nu(x)$ of $\tau$ with minimal $i$-value is as claimed. To see this, let $A$ be an arbitrary apartment containing $x$ and $z_0$. Then $\tau$ is contained in $A$ because this is true for $[z_0,x]$. Explicitly it can be described as follows. Choose an identification of $A$ with $\mathbb{R}^{d+1}/{\mathbb{R}}.(1,\ldots,1)$ as before, with $z_0\in A$ corresponding to the class of the origin in $\mathbb{R}^{d+1}/{\mathbb{R}}.(1,\ldots,1)$. After reindexing the basis if necessary there are sequences $0\le r_0<r_1<\ldots<r_s=d$ and $0<m_1<\ldots<m_s$ (some $1\le s\le d$) such that the vertex $x$ is represented by $y_s=\sum_{i=1}^s\sum_{j=r_{i-1}+1}^{r_i}m_ie_{{j}}$. Then $\{hy_s|\,\,0\le h\le 1\}$ represents $[z_0,x]$, and $\tau$ is a $s$-simplex, the other vertices are represented by $y_t=\sum_{i=1}^t\sum_{j=r_{i-1}+1}^{r_i}m_ie_{{j}}+\sum_{i=t+1}^s\sum_{j=r_{i-1}+1}^{r_i}(m_i-1)e_{{j}}$ for $t=0,\ldots, s-1$. In particular $\nu(x)$ is represented by $y_0$ and it is clear that it has minimal $i$-value among the vertices of $A$ incident to $x$. Since any vertex in $X$ incident to $x$ lies in such an $A$ the assertion follows. Also the other claims can immediately be read off from this analysis on an apartment.

\begin{pro}\label{baswofu} Let $\widehat{\eta}=(x_1,\ldots,x_k)$ be as in \ref{iordpoin} and let $\eta=\{x_1,\ldots,x_k\}$. Then$$M_0=\{u\in X^0|\,\,\{u\}\cup\eta\,\mbox{ is a }\,\,k\mbox{-simplex and}\,\,i(u)<i(x_1)\}$$is contained in $N_{\widehat{\eta}}$ and stable with respect to $\widehat{\eta}$.\end{pro}

{\sc Proof:} The containment $M_0\subset N_{\widehat{\eta}}$ follows from \ref{iordpoin} with $k$ instead of $(k-1)$. To prove that $M_0$ is stable with respect to $\widehat{\eta}$ let $u_1, u_2\in M_0$. First it follows from \ref{iordpoin} (with $k$ and $k+1$ instead of $(k-1)$) and \ref{minnei} that $\nu(x_k)\in W^{\widehat{\eta}}_{u_1,u_2}$, hence $[\widehat{\eta}|u_1,u_2]$ is defined. Since $([\widehat{\eta}|u_1,u_2],x_1,\ldots,x_k)$ is a pointed $k$-simplex and since for any pointed $k$-simplex the underlying {\it cyclic} ordering of the vertices is independent of the pointing, there are --- in view of \ref{iordpoin} with $k$ instead of $(k-1)$ --- only the two possibilities $i([\widehat{\eta}|u_1, u_2])<i(x_1)$ and $i([\widehat{\eta}|u_1, u_2])>i(x_k)$. If $i([\widehat{\eta}|u_1, u_2])<i(x_1)$ then $[\widehat{\eta}|u_1, u_2]\in M_0$ and we are done. If $i([\widehat{\eta}|u_1, u_2])>i(x_k)$ then $\nu([\widehat{\eta}|u_1, u_2])\in W^{\widehat{\eta}}_{u_1,u_2}$ by \ref{iordpoin} and \ref{minnei}. Moreover $$\ell((\nu([\widehat{\eta}|u_1, u_2]),\widehat{\eta}))\le \ell(([\widehat{\eta}|u_1, u_2],\widehat{\eta}))$$also by \ref{minnei}. Since $\nu([\widehat{\eta}|u_1, u_2])\ne[\widehat{\eta}|u_1, u_2]$ this contradicts the definition of $[\widehat{\eta}|u_1, u_2]$. Hence $i([\widehat{\eta}|u_1, u_2])>i(x_k)$ can not happen and the proof is finished.\\

{\sc Proof of Theorem \ref{acykrit}:} Given $\eta, \eta'\in X^{k-1}$ let $x_1,\ldots,x_k$, resp. $x'_1,\ldots,x'_k$, be that enumeration of the vertices of $\eta$, resp. of $\eta'$, which satisfies $i(x_1)<\ldots<i(x_k)$, resp. $i(x'_1)<\ldots<i(x'_k)$. We define $$\eta\widetilde{<}\eta'\Longleftrightarrow\left\{\begin{array}{c}\mbox{there is a}\,\,1\le q\le k\,\,\mbox{such that}\\i(x_t)=i(x'_t)\,\,\mbox{for}\,\,1\le t\le q-1\,\,\mbox{and}\,\,i(x_q)<i(x'_q)\end{array}\right\}.$$We define $\eta\widetilde{=}\eta'$ if $i(x_t)=i(x'_t)$ for all $1\le t\le k$. Let ${\nabla}:X^{k-1}\longrightarrow\mathbb{N}$ be the surjective map with $$\nabla(\eta)<\nabla(\eta')\Longleftrightarrow\eta\widetilde{<}\eta'$$$$\nabla(\eta)=\nabla(\eta')\Longleftrightarrow \eta\widetilde{=}\eta'.$$For a $k$-simplex $\sigma\in X^k$ let $\sigma^-\in X^{k-1}$ be the $(k-1)$-simplex obtained from $\sigma$ by omitting the vertex $x\in\sigma$ for which $i(x)$ is minimal. We need to show that$$\prod_{\eta\in X^{k-1}}{\mathcal F}(\eta)\stackrel{\partial^k}{\longrightarrow}\prod_{\sigma\in X^k}{\mathcal F}(\sigma)\stackrel{\partial^{k+1}}{\longrightarrow}\prod_{\tau\in X^{k+1}}{\mathcal F}(\tau)$$is exact. So let a $k$-cocycle $c=(c_{\sigma})_{\sigma\in X^k}\in\ke(\partial^{k+1})$ be given. It suffices to show that there is a sequence of $(k-1)$-cochains $(b_n)_{n\in\mathbb{N}}=((b_{n,\eta})_{\eta\in X^{k-1}})_{n\in\mathbb{N}}$ with $b_{n,\eta}\in{\mathcal F}(\eta)$ satisfying the following properties:\begin{description}\item[(i)] $b_{n,\eta}=b_{\nabla(\eta),\eta}$ for all $n\ge \nabla(\eta)$.\item[(ii)] $b_{n,\eta}=0$ for all $\eta\in X^{k-1}$ with $\nabla(\eta)>n$.\item[(iii)] For all $\sigma\in X^{k}$ with $\nabla(\sigma^-)\le n$ we have $(\partial^kb_n-c)_{\sigma}=0$.\end{description}Then the cochain $b_{\infty}=(b_{\infty,\eta})_{\eta\in X^{k-1}}$ defined by $b_{\infty,\eta}=b_{\nabla(\eta),\eta}$ will be a preimage of $c$, as follows from (i) and (iii).

We construct $(b_n)_{n\in\mathbb{N}}$ inductively. Suppose $b_{n-1}$ has been constructed. We set $b_{n,\eta}=b_{n-1,\eta}$ for all $\eta\in X^{k-1}$ with $\nabla(\eta)<n$, and $b_{n,\eta}=0$ for all $\eta\in X^{k-1}$ with $\nabla(\eta)>n$. Now suppose we have $\eta\in X^{k-1}$ with $\nabla(\eta)=n$. Let $x_1,\ldots,x_k$ be that enumeration of the vertices of ${\eta}$ which satisfies $i(x_1)<\ldots<i(x_k)$. Consider the set \begin{gather}M_0=\{z\in X^0|\,\,\{z\}\cup\eta\,\mbox{ is a }\,k\mbox{-simplex and}\,i(z)<i(x_1)\}.\label{petdis}\end{gather}We know from \ref{iordpoin} and \ref{baswofu} that $\widehat{\eta}=(x_1,\ldots,x_k)$ is a pointed $(k-1)$-simplex and that $M_0$ is stable with respect to $\widehat{\eta}$ and contained in $N_{\widehat{\eta}}$. If $M_0=\emptyset$ we put $b_{n,\eta}=0$. So assume now that $M_0\ne\emptyset$. Let $z$, $z'\in M_0$ with $\{z,z'\}\in X^1$. We compute (with $\pm$, resp. $r$, denoting the respective incidence numbers, resp. restriction maps):\begin{align} \quad{\,\,\,\,\,\,} & \partial^{k+1}(((\partial^kb_{n-1}-c)_{z''\cup\eta})_{z''\in M_0})_{\{z,z'\}\cup\eta\}}\notag \\ = &\pm\,\, r((\partial^kb_{n-1}-c)_{z\cup\eta})+ \pm\,\, r((\partial^kb_{n-1}-c)_{z'\cup\eta})\notag \\= & \pm\,\, r((\partial^kb_{n-1}-c)_{z\cup\eta})+ \pm\,\, r((\partial^kb_{n-1}-c)_{z'\cup\eta})+\sum_{\{z,z'\}\subset\sigma\subset \{z,z'\}\cup\eta}\pm\,\, r((\partial^kb_{n-1}-c)_{\sigma})\notag \\ = &  \sum_{\sigma\subset \{z,z'\}\cup\eta}\pm\,\, r((\partial^kb_{n-1}-c)_{\sigma})\notag \\ = & (\partial^{k+1}(\partial^kb_{n-1}-c))_{\{z,z'\}\cup\eta}\notag\end{align}and this is zero because $c$ is a cocycle. For the second equality note that for all $\sigma\in X^k$ with $\{z,z'\}\subset\sigma\subset \{z,z'\}\cup\eta$ we have $\nabla(\sigma^-)<n$ which by induction hypothesis implies $(\partial^k b_{n-1}-c)_{\sigma}=0$. We have seen that $((\partial^k b_{n-1}-c)_{\{z\}\cup\eta})_{z\in M_0}$ lies in$$\ke[\prod_{z\in M_0}{\mathcal 
F}(\{z\}\cup\eta)\stackrel{\partial^{k+1}}{\longrightarrow}\prod_{z, z'\in M_0\atop \{z,z'\}\in X^1}{\mathcal 
F}(\{z,z'\}\cup\eta)].$$We can therefore define $b_{n,\eta}\in {\mathcal 
F}(\eta)$ as a preimage of $[\sigma:\sigma^-]((c-\partial^k b_{n-1})_{\{z\}\cup\eta})_{z\in M_0}$, by hypothesis ${\mathcal S}(k)$. To see that $b_n$ satisfies (iii) for $\sigma\in X^k$ with $\nabla(\sigma^-)=n$ we compute\begin{align}(\partial b_n-c)_{\sigma} & = (\sum_{\eta\subset\sigma}[\sigma:\eta]r_{\sigma}^{\eta}(b_{n,\eta}))-c_{\sigma}\notag \\
{} & = (\sum_{\eta\subset\sigma\atop \nabla(\eta)<n}[\sigma:\eta]r_{\sigma}^{\eta}(b_{n,\eta}))+[\sigma:\sigma^-]r_{\sigma}^{\sigma^-}(b_{n,\sigma^-})-c_{\sigma}\notag \\
{} & = (\sum_{\eta\subset\sigma\atop \nabla(\eta)<n}[\sigma:\eta]r_{\sigma}^{\eta}(b_{n-1,\eta}))+(c-\partial^k b_{n-1})_{\sigma}-c_{\sigma}\notag \\
{} & = (\sum_{\eta\subset\sigma\atop \nabla(\eta)<n}[\sigma:\eta]r_{\sigma}^{\eta}(b_{n-1,\eta}))- (\sum_{\eta\subset\sigma}[\sigma:\eta]r_{\sigma}^{\eta}(b_{n-1,\eta}))
\notag\end{align}and this is zero because we have $b_{n-1,\sigma^-}=0$ by induction hypothesis (ii).\\ 

The reader will have observed that any $\mathbb{N}$-valued function $i$ on $X^0$ which takes different values on incident vertices gives rise to a local vanishing criterion like ${\mathcal{S}}(k)$, by the same formal proof above. However, the applicability of the resulting criterion depends on the control one gets over the corresponding sets $M_0$ defined analogously through formula (\ref{petdis}). In this optic, the virtue of our particular choice of $i$ lies in the fact that we can control the corresponding sets $M_0$: they satisfy the {\it local} (no reference to the global function $i$) property of being stable with respect to $\widehat{\eta}$; hence our vanishing criterion ${\mathcal{S}}(k)$, expressed entirely in local terms.\\

A {\it homological coefficient system} (hcs) ${\mathcal V}$ of abelian groups on a building $X$ is the assignment of an abelian group ${\mathcal V}(\tau)$ to every simplex ${\tau}$ of $X$, and a homomorphism $r^{{\tau}}_{{\sigma}}:{\mathcal V}({\tau})\to{\mathcal V}({\sigma})$ to every face inclusion ${\sigma}\subset{\tau}$, such that $r_{{\rho}}^{{\tau}}\circ r_{{\tau}}^{{\sigma}}=r_{{\rho}}^{{\sigma}}$ whenever ${\rho}\subset{\tau}\subset{\sigma}$, and $r_{{\tau}}^{{\tau}}$ is the identity.

The group $C_k(X,{\mathcal V})$ of $k$-chains ($0\le k\le d$) consists of all the finitely supported maps $c$ which assign to each $k$-simplex ${\tau}$ an element $c_{{\tau}}\in{\mathcal V}({\tau})$. Define$$\partial=\partial_{k}:C_{k+1}(X,{\mathcal V})\longrightarrow C_{k}(X,{\mathcal V})$$by the rule$$(\partial c)_{{\tau}}=\sum_{{\tau'}\supset{\tau}}[\tau':\tau]r_{{\tau}}^{{\tau'}}(c_{{\tau'}}).$$Then $(C_{\bullet}(X,{\mathcal V}),\partial)$ is a complex ($\partial^2=0$), and its homology groups are denoted $H_k(X,{\mathcal V})$.\\

Consider for $1\le k\le d$ the following condition ${\mathcal S}^*(k)$ for a hcs ${\mathcal V}$ on $X$: For any pointed $(k-1)$-simplex $\widehat{\eta}$ with underlying $(k-1)$-simplex $\eta$ and for any subset $M_0$ of $N_{\widehat{\eta}}$ which is stable with respect to $\widehat{\eta}$, the following subquotient complex of $(C_{\bullet}(X,{\mathcal V}),\partial)$ is exact:$$\bigoplus_{z,z'\in {M_0}\atop\{z,z'\}\in{X^1}}{\mathcal 
V}(\{z,z'\}\cup\eta)\stackrel{\partial_{k}}{\longrightarrow}\bigoplus_{z\in {M}_0}{\mathcal 
V}(\{z\}\cup\eta)\stackrel{\partial_{k-1}}{\longrightarrow}{\mathcal V}(\eta).$$

\begin{satz}\label{hoacykrit} Let ${\mathcal V}$ be a hcs on $X$. Let $1\le k\le d$ and suppose ${\mathcal S}^*(k)$ holds true. Then $H_k(X,{\mathcal V})=0$.  
\end{satz}

{\sc Proof:} We need to show that$$\bigoplus_{\tau\in X^{k+1}}{\mathcal V}(\tau)\stackrel{\partial_k}{\longrightarrow}\bigoplus_{\sigma\in X^k}{\mathcal V}(\sigma)\stackrel{\partial_{k-1}}{\longrightarrow}\bigoplus_{\eta\in X^{k-1}}{\mathcal V}(\eta)$$is exact. We use notations from the proof of \ref{acykrit}. For $n\in \mathbb{Z}_{\ge0}$ and elements $c=(c_{\sigma})_{\sigma}\in \oplus_{\sigma\in X^k}{\mathcal V}(\sigma)$ consider the condition $$C(n):\,\,\mbox{for all}\,\,\sigma\in X^k\,\,\mbox{with}\,\,\nabla(\sigma^-)\ge n\,\,\mbox{we have}\,\,c_{\sigma}=0.$$Similarly as in the proof of \ref{acykrit} one shows by induction on $n$: all elements $c\in\ke(\partial_{k-1})$ which satisfy $C(n)$ ly in $\bi(\partial_k)$. Indeed, given such an element $c$ one uses ${\mathcal S}^*(k)$ in order to modify $c$ by an element of $\bi(\partial_k)$ in such a way that it even satisfies $C(n-1)$, and then the induction hypothesis applies.  

\section{$p$-adic hyperplane arrangements}

\label{sechyar}

Let $K$ denote a non-archimedean locally compact field, ${\cal O}_K$ its ring of integers, $\pi\in{\cal O}_K$ a fixed prime element and $k$ the residue field. Let $X$ be the Bruhat-Tits building of ${\rm PGL}_{d+1}/K$; it has Coxeter data of type $\widetilde{A}_d$. A concrete descriptions of $X$ is the following. A lattice in the $K$-vector space $K^{d+1}$ is a free ${\cal O}_K$-submodule of $K^{d+1}$ of rank $d+1$. Two lattices $L$, $L'$ are called homothetic if $L'=\lambda L$ for some $\lambda\in K^{\times}$. We denote the homothety class of $L$ by $[L]$. The set of vertices of $X$ is the set of the homothety classes of lattices (always: in $K^{d+1}$). For a lattice chain$$\pi L_k\subsetneq L_{0}\subsetneq L_{1}\subsetneq\ldots\subsetneq L_k$$we declare $([L_0],\ldots,[L_k])$ to be a pointed $k$-simplex. This defines a simplicial structure with orientation.

The following definitions are due to de Shalit \cite{ds} (who takes $R=K$, $\kara(K)=0$ below). Let ${\mathcal{A}}$ be a non empty subset of $\mathbb{P}(K^{d+1})$. We view ${\mathcal{A}}$ as a set of lines in $K^{d+1}$, or hyperplanes in $(K^{d+1})^*$. We write $e_a$ for the line represented by $a\in K^{d+1}-\{0\}$, so that $e_{\lambda a}=e_a$ for any $\lambda\in K^{\times}$.

Let $R$ be a commutative ring and let $\widetilde{E}$ be the free exterior algebra over $R$, on the set ${\mathcal{A}}$. It is graded and anti-commutative. There is a unique derivation $\delta:\widetilde{E}\to \widetilde{E}$, homogeneous of degree $-1$, mapping each $e\in{\mathcal{A}}$ to $1$. It satisfies $\delta^2=0$ and \begin{gather}\delta(e_{0}\wedge\ldots\wedge e_{k})=\sum_{i=0}^k(-1)^ie_{0}\wedge\ldots\wedge\widehat{e_{i}}\wedge\ldots\wedge e_{k}\notag\\\bi(\delta)=\ke(\delta).\notag\end{gather}The subalgebra $E=\bi(\delta)=\ke(\delta)$ of $\widetilde{E}$ is generated by all elements $e-e'$ for $e,e'\in{\mathcal{A}}$. There is an exact sequence\begin{gather}0\longrightarrow E\longrightarrow \widetilde{E}\stackrel{\delta}{\longrightarrow} E[1]\longrightarrow 0\label{espli}\end{gather}and any $e\in{\mathcal{A}}$ supplies a splitting $E[1]\to\widetilde{E}$, $x\mapsto e\wedge x$. 

Let $x\in X^0$ be a vertex. We say that an element $g\in\widetilde{E}$ is a {\it standard generator of} $I(x)$ if there are a lattice $L_x$ representing $x$ and elements $\{a_0,\ldots,a_m\}$ of ${\mathcal A}\cap L_x-\pi L_x$, linearly dependent modulo $\pi L_{x}$, such that $g=\delta(e_{a_0}\wedge\ldots\wedge e_{a_m})$. We define the ideal $I(x)$ in $\widetilde{E}$ as the one generated by the standard generators of $I(x)$. For an arbitrary simplex $\sigma$ we define the ideal $I(\sigma)=\sum_{x\in\sigma}I(x)$. We set$$\widetilde{A}(\sigma)=\frac{\widetilde{E}}{I(\sigma)},\quad\quad\quad{A}(\sigma)=\frac{{E}}{E\cap I(\sigma)}.$$The split exact sequence (\ref{espli}) provides us with a split exact sequence\begin{gather}0\longrightarrow A(\sigma)\longrightarrow \widetilde{A}(\sigma)\stackrel{\delta}{\longrightarrow}A(\sigma)[1]\longrightarrow 0.\label{aspli}\end{gather}The ideal $I(\sigma)$ is homogeneous, hence there is a natural grading on $\widetilde{A}(\sigma)$ and on $A(\sigma)$. We denote by $\widetilde{A}^q(\sigma)$ resp. by ${A}^q(\sigma)$ the $q$-th graded piece. For varying $\sigma$ the $\widetilde{A}(\sigma)$ and ${A}(\sigma)$ form ccs $\widetilde{A}$ and $A$ on $X$.\\

Suppose we are given a lattice chain$$\pi L_k\subsetneq L_{1}\subsetneq L_{2}\subsetneq\ldots\subsetneq L_k.$$Let $x_j\in{X^0}$ be the vertex defined by $L_j$; then $\widehat{\eta}=(x_1,\ldots,x_k)$ is a pointed $(k-1)$-simplex; we denote by $\eta$ the underlying non pointed $(k-1)$-simplex. We write $L_{0}=\pi L_k$. As long as $L_k$ is fixed we abuse notation in that we identify an element $e\in{\mathcal{A}}$ with an element $a\in L_k-L_0$ for which $e=e_a$; such an $a$ is unique up to a unit in ${\mathcal O}_K$. Thus we regard ${\mathcal {A}}$ as a subset of $L_k-L_0$.

Assume that ${\mathcal{A}}$ is finite and fix a linear ordering $\prec$ on ${\mathcal{A}}$ which is {\it adapted to $\widehat{\eta}$}. By definition, this means $\max({\mathcal{A}}\cap L_j-L_{j-1})\prec\min({\mathcal{A}}\cap L_{j+1}-L_{j})$ for all $1\le j\le k-1$. For $S\subset{\mathcal{A}}$ and $e\in S$ we define the ${\mathcal O}_K$-submodule $L(\widehat{\eta},\prec,S,e)$ of $K^{d+1}$ as follows: Let $1\le j\le k$ be the number for which $e\in L_j-L_{j-1}$; then $$L(\widehat{\eta},\prec,S,e)=<L_{j-1},\,\,\{e'\in S|\,\,e'\prec e\,\,\mbox{or}\,\,e'=e\}>_{{\mathcal O}_K}.$$We say that {\it $e$ is $(S,\widehat{\eta})$-special with respect to $\prec$} if$$e=\max_{\prec}({\mathcal{A}}\cap L(\widehat{\eta},\prec,S,e))$$(Here $\max_{\prec}(Q)$ for a subset $Q$ of ${\mathcal{A}}$ means the maximal element of $Q$ with respect to the fixed ordering $\prec$. The subscript $\prec$ does {\it not} indicate that $\prec$ is a running parameter.)

Fix another linear ordering $<$ on ${\mathcal{A}}$ and for subsets $S$ of ${\mathcal{A}}$ let $e_S=e_{0}\wedge\ldots\wedge e_{r}\in \widehat{E}$ where $e_0<e_1<\ldots<e_r$ is the increasing enumeration of the elements of $S$ in the ordering $<$. (The ordering $<$ may be taken to be $\prec$, but the role of these two orderings will be completely unrelated in the following).

\begin{lem}\label{broci} $\widetilde{A}(\eta)$ is a free $R$-module of finite rank, a basis is the set$$\{e_S|\,\,S\subset {\mathcal{A}},\,\,\mbox{all}\,\, e\in S\,\,\mbox{are}\,\,(S,\widehat{\eta})\mbox{-special with respect to} \prec\}.$$
\end{lem}

{\sc Proof:} This is \cite{ds} Theorem 2.5: the "broken circuit theorem"\, (there $R=K$, $\kara(K)=0$ and $\prec=<$, but this is irrelevant).\\

Let $N_{\widehat{\eta}}$ be as in section \ref{thec}. A subset $M_0\subset N_{\widehat{\eta}}$ corresponds to a collection of lattices $(L_z)_{z\in M_0}$ with $L_{0}\subsetneq L_z\subsetneq L_1$ for all $z\in M_0$. The following lemma is clear. 

\begin{lem}\label{fulat} $M_0$ is stable with respect to $\widehat{\eta}$ if and only if for all $z_1, z_2\in M_0$ the lattice $L_{z_1}\cap L_{z_2}$ represents an element of $M_0$.
\end{lem}

Suppose that $M_0\subset N_{\widehat{\eta}}$ is stable with respect to $\widehat{\eta}$. In particular there is a $x_0\in M_0$ with $L_{x_0}\subset L_z$ for all $z\in M_0$. We say that a collection $(\prec_z)_{z\in{M_0}}$, indexed by $M_0$, of linear orderings on ${\mathcal{A}}$ is {\it adapted to $(M_0,\widehat{\eta})$} if the following conditions hold:\begin{itemize}\item For any $z\in M_0$ the ordering $\prec_z$ is adapted to the pointed $k$-simplex $(z,\widehat{\eta})$.\item For any $z_1, z_2\in M_0$ with $L_{z_1}\subset L_{z_2}$ we have $[e\prec_{z_1}e'\Leftrightarrow e\prec_{z_2}e']$ for all pairs $e,e'\in{\mathcal{A}}\cap L_{z_1}$ and also for all pairs $e,{e'}\in{\mathcal{A}}\cap L_{z_2}-L_{z_1}$.\item For any $z\in M_0$ we have $[e\prec_{x_0}e'\Leftrightarrow e\prec_{z}e']$ for all pairs $e,e'\in{\mathcal{A}}\cap L_{k}-L_{z}$.\item We have $e\prec_{x_0}e'$ for all pairs $e,e'\in{\mathcal{A}}$ with $e\in\bigcup_{z\in M_0}L_z$ and ${e'}\notin\bigcup_{z\in M_0}L_z$.\end{itemize}

\begin{lem}\label{collex} Collections of linear orderings on ${\mathcal{A}}$ adapted to $(M_0,\widehat{\eta})$ exist.
\end{lem}

{\sc Proof:} The referee suggested the following proof (our original one was unnecessarily complicated). Let $U=\bigcup_{z\in M_0}L_z$. Fix a linear ordering $\prec$ on ${\mathcal{A}}$ adapted to $\widehat{\eta}$ which satisfies $e\prec e'$ for all pairs $e,e'\in{\mathcal{A}}$ with $e\in U$ and ${e'}\notin U$. For a $z\in M_0$ let $\prec_z$ be the ordering which satisfies firstly $e\prec_z e'\prec_z e''$ for all triples $e,e',e''\in{\mathcal{A}}$ with $e\in L_z$, with $e'\in U-L_z$ and with $e''\in L_k-U$, and secondly $[e\prec e'\Leftrightarrow e\prec_{z}e']$ for all pairs $e,e'\in{\mathcal{A}}\cap L_{z}$, for all pairs $e,e'\in{\mathcal{A}}\cap U-L_{z}$ and for all pairs $e,e'\in{\mathcal{A}}-({\mathcal{A}}\cap U)$. It is straightforwardly checked that $(\prec_z)_{z\in{M_0}}$ is adapted to $(M_0,\widehat{\eta})$.\\

Fix a collection $(\prec_z)_{z\in{M_0}}$ of linear orderings on ${\mathcal{A}}$ adapted to $(M_0,\widehat{\eta})$. Let $$G(\widehat{\eta};M_0)=\{\,\,e_S\,\,|\begin{array}{c}\,\,S\subset{\mathcal{A}},\,\,\mbox{for all}\,\,e\in S\,\,\mbox{there is a}\\\,\,z\in M_0\,\,\mbox{such that}\,\,e\,\,\mbox{is}\,\,(S,(z,\widehat{\eta}))\mbox{-special}\end{array}\}$$where $(S,(z,\widehat{\eta}))$-speciality is to be understood with respect to $\prec_z$. Let $J(\widehat{\eta};M_0)\subset \widetilde{E}$ be the ideal generated by the set$$I(\eta)\,\,\bigcup\,\,\{g\in\widetilde{E}|\,\,g\,\,\mbox{is a standard generator of}\,\,I(z)\,\,\mbox{for all}\,\,z\in M_0\}.$$

\begin{pro}\label{erzgrad} $G(\widehat{\eta};M_0)$ is finite and generates $\widetilde{E}/J(\widehat{\eta};M_0)$ as an $R$-module.
\end{pro}

{\sc Proof:} Clearly the set of {\it all} $e_S$ with $S\subset {\mathcal{A}}$ is generating. Now suppose that $S\subset {\mathcal{A}}$ does not satisfy the condition defining $G(\widehat{\eta};M_0)$. That is, there exists an $e\in S$ such that for all $z\in M_0$ this $e$ is {\it not} $(S,(z,\widehat{\eta}))$-special. Fix such an $e$.

We first claim that there is an $\widehat{e}\in S$ and a subset $\widehat{S}\subset S$ such that for all $z\in M_0$ the following statements (1) and (2) are satisfied:\begin{description}\item[(1)] For all $e'\in \widehat{S}$ we have $e'\prec_z e\prec_z\widehat{e}$.\item[(2)] The element $\delta(e_{\widehat{S}\cup\{e,\widehat{e}\}})$ belongs to $J(\widehat{\eta};M_0)$.\end{description}

We distinguish three cases. First consider the case $e\in L_j-L_{j-1}$ for some $2\le j\le k$. Then we set $\widehat{S}=\{e'\in S|\,\,e'\prec_{x_0}e\,\,\mbox{and}\,\,e'\notin L_{j-1}\}$ and $$\widehat{e}=\max_{\prec_{x_0}}({\mathcal{A}}\cap L((x_0,\widehat{\eta}),\prec_{x_0},S,e)).$$ The fact that $e$ is not $(S,(x_0,\widehat{\eta}))$-special and the properties of the adapted collection $(\prec_z)_{z\in{M_0}}$ give statement (1). Moreover $\widehat{S}\cup\{e,\widehat{e}\}$ is linearly dependent modulo $L_{j-1}$, hence $\delta(e_{\widehat{S}\cup\{e,\widehat{e}\}})$ lies in $I(x_{j-1})\subset I(\eta)$ and we get statement (2).

Now consider the case $e\in L_1-\bigcup_{z\in M_0}L_z$. Then we set $\widehat{S}=\{e'\in S|\,\,e'\prec_{x_0}e\}$ and $\widehat{e}=\max_{\prec_{x_0}}({\mathcal{A}}\cap L((x_0,\widehat{\eta}),\prec_{x_0},S,e))$. Again the fact that $e$ is not $(S,(x_0,\widehat{\eta}))$-special and the properties of the adapted collection $(\prec_z)_{z\in{M_0}}$ give statement (1). For $z\in M_0$ let $\widehat{S}_{z}=\widehat{S}-(\widehat{S}\cap L_{z})$. The fact that $e$ is not $(S,(z,\widehat{\eta}))$-special tells us that $\widehat{S}_z\cup\{e,\widehat{e}\}$ is linearly dependent modulo $L_{z}$. But then also the subset of $\pi^{-1}L_z-L_z$ which defines the same elements in ${\mathbb P}(K^{d+1})$ as does $\widehat{S}\cup\{e,\widehat{e}\}$ is linearly dependent modulo $L_{z}$, because it contains $\widehat{S}_z\cup\{e,\widehat{e}\}$. Therefore $\delta(e_{\widehat{S}\cup\{e,\widehat{e}\}})$ is a standard generator of $I(z)$. We have shown statement (2).

Finally the case $e\in\bigcup_{z\in M_0}L_z$. Let $z_e\in M_0$ be such that $e\in L_{z_e}$ and $L_{z_e}$ is minimal with this property (this $z_e$ is unique since $M_0$ is stable). We set $\widehat{S}=\{e'\in S|\,\,e'\prec_{z_e}e\}$ and $\widehat{e}=\max_{\prec_{z_e}}({\mathcal{A}}\cap L((z_e,\widehat{\eta}),\prec_{z_e},S,e))$. The fact that $e$ is not $(S,(z_e,\widehat{\eta}))$-special gives statements (1) and (2) in this case (here $\delta(e_{\widehat{S}\cup\{e,\widehat{e}\}})\in I(x_{k})\subset I(\eta)$). Again this is straightforwardly checked using the properties of the adapted collection $(\prec_z)_{z\in{M_0}}$.

The claim established, statement (2) tells us that we may replace $e_{\widehat{S}\cup\{e\}}$ in $\widetilde{E}/J(\widehat{\eta};M_0)$ by a linear combination of elements $e_{S'}$ with each $S'$ arising from $\widehat{S}\cup\{e,\widehat{e}\}$ by deleting an element of $\widehat{S}\cup\{e\}$. By statement (1) of our claim this means that we may replace $e_S$ by a linear combination of elements $e_{S''}$ with each $S''$ satisfying $S\prec_z S''$ for all $z\in M_0$ (for the lexicographic ordering $\prec_z$ on the set of subsets of ${\mathcal{A}}$ of fixed cardinality derived from the ordering $\prec_z$ on ${\mathcal{A}}$). Repeating the process proves that $G(\widehat{\eta};M_0)$ is generating. That it is finite follows from \ref{broci}. We are done.\\

We define the complex$$K(\widehat{\eta},M_0)=\left\{\begin{array}{c}0\longrightarrow\frac{\widetilde{E}}{\bigcap_{z\in M_0}I(\eta\cup\{z\})}\longrightarrow\\\prod_{z\in M_0}\widetilde{A}(\eta\cup\{z\})\longrightarrow\prod_{z_1,z_2\in M_0\atop \{z_1,z_2\}\in{X^1}}\widetilde{A}(\eta\cup\{z_1,z_2\})\longrightarrow\ldots\end{array}\right\}.$$

\begin{pro}\label{simind} For non-empty $M_0$ as above the following statements hold:\\ ${\bf (a)}$ The complex $K(\widehat{\eta},M_0)$ is exact.\\${\bf (b)}$$$J(\widehat{\eta};M_0)=\bigcap_{z\in M_0}I(\eta\cup\{z\}).$$\\${\bf (c)}$ $\widetilde{E}/\bigcap_{z\in M_0}I(\eta\cup\{z\})$ is a free $R$-module and $G(\widehat{\eta};M_0)$ is a basis.
\end{pro}

{\sc Proof:} For $n\ge1$ let ${\bf (a)_n}$, ${\bf (b)_n}$ and ${\bf (c)_n}$ be the corresponding statements for all $M_0$ with $1\le|M_0|\le n$. We will prove these statements by simultanuous induction on $n$. 
 
Statements ${\bf (a)_1}$ and ${\bf (b)_1}$ are clear, and ${\bf (c)_1}$ is \ref{broci}. Now we assume $n>1$. First we prove ${\bf (a)_n}$. Choose a $y\in M_0$ for which $L_y$ is maximal, i.e. such that there is no $z\in M_0$ with $L_y\subsetneq L_z$. If we set $$M_1=M_0-\{y\},\quad\quad M_1'=\{z\in M_1|\quad L_z\subset L_y\}$$then $M_1$ and $M_1'$ are stable with respect to $\widehat{\eta}$ resp. $(y,\widehat{\eta})$. We claim\begin{gather}\bigcap_{z'\in M_1'}I(\eta\cup\{y,z'\})\,\,\subset\,\, I(\eta\cup\{y\})+\bigcap_{{z}\in M_1}I(z).\label{keycla}\end{gather}For $x\in M_0$ and $S\subset{\mathcal{A}}$ let $S^x$ be the uniquely determined subset of $L_x-\pi L_x$ which determines the same set of lines in $K^{d+1}$ (or: the same subset of $\mathbb{P}(K^{d+1}))$ as does $S$ (recall that by convention we view $S$ as a subset of $L_k-L_0$). By induction hypothesis ${\bf (b)_{n-1}}$, applied to the pointed $k$-cell $(y,\widehat{\eta})$ and $M_1'\subset N_{(y,\widehat{\eta})}$, a typical generator $g$ of the left hand side of (\ref{keycla}) satisfies at least one of the following conditions:\begin{description}\item[(i)] $g\in I(\eta\cup\{y\})$\item[(ii)] $g=\delta(e_{0}\wedge\ldots\wedge e_{r})$ for some $S=\{e_0,\ldots,e_r\}\subset{\mathcal{A}}$ such that $S^{z'}$ is linearly dependent modulo $\pi L_{z'}$ for all $z'\in M_1'$.\end{description}

To show that $g$ lies in the right hand side of (\ref{keycla}) is difficult only if (i) does not hold. In that case we claim $g\in\bigcap_{z\in M_1}I(z)$. So let $S$ be as in (ii) and let $z\in M_1$ be given. Since $M_0$ is stable with respect to $\widehat{\eta}$ there exists a $z'\in M_1'$ such that $L_y\cap L_z=L_{z'}$. Now $S^{z'}$ is linearly dependent modulo $\pi L_{z'}$. If $S^y$ was linearly dependent modulo $\pi L_y$ we would have $g\in I(y)$, but then (i) would hold. Hence $S^y$ is linearly independent modulo $\pi L_y$, hence so is $S^{z'}\cap S^y$, and hence $\widetilde{S}=S^{z'}-(S^{z'}\cap S^y)$ is linearly dependent modulo $\pi L_z$ (since $S^{z'}$ is). On the other hand $\widetilde{S}\subset \pi L_y$ and this implies $\widetilde{S}\cap \pi L_z=\emptyset$ (since $\pi L_y\cap \pi L_z=\pi L_{z'}$). Thus $\widetilde{S}$ is a subset of $L_z-\pi L_z$ which is linearly dependent modulo $\pi L_{z}$. In particular, $S^z$ is linearly dependent modulo $\pi L_z$ (since $\widetilde{S}\subset S^z$), hence $g\in I(z)$.

We have established the claim (\ref{keycla}). Next we claim that\begin{gather}\widetilde{E}\longrightarrow\prod_{z\in {M}_0}\widetilde{A}(\eta\cup\{z\})\longrightarrow\prod_{z_1,z_2\in {M_0}\atop\{z_1,z_2\}\in{X^1}}\widetilde{A}(\eta\cup\{z_1,z_2\})\label{zielsu}\end{gather}is exact; this is equivalent with exactness of $K(\widehat{\eta},M_0)$ at the first non trivial degree. Let $(s_z)_{z\in M_0}$ be an element of the kernel of the second arrow in (\ref{zielsu}). By induction hypothesis ${\bf (a)_{n-1}}$ for $M_1$ we may assume, after modifying $(s_z)_{z\in M_0}$ by the image of an element of $\widetilde{E}$, that $s_z=0$ for all $z\in M_1$. Then it follows that $$s_{y}\in\ke[\widetilde{A}(\eta\cup\{y\})\longrightarrow\prod_{{z'\in {M_1'}}}\widetilde{A}(\eta\cup\{z',y\})].$$This means that $s_y$ can be lifted to an element $\tilde{s}_y$ of the left hand side of (\ref{keycla}). By (\ref{keycla}) there exist $\tilde{b}\in I(\eta\cup\{y\})$ and $\tilde{c}\in\bigcap_{{z}\in M_1}I(\eta\cup\{z\})$ with $\tilde{s}_y=\tilde{b}+\tilde{c}$. Thus $\tilde{c}$ is a preimage of $(s_z)_{z\in M_0}$ and the exactness of (\ref{zielsu}) is proven. If we define the complex$$K^y(\widehat{\eta},M_0)=\left\{\begin{array}{c}0\longrightarrow\frac{\widetilde{E}}{\bigcap_{z\in M_1}I(\eta\cup\{z\})}\longrightarrow\\\frac{\widetilde{E}}{\bigcap_{z'\in M_1'}I(\eta\cup\{y,z'\})}\times\prod_{z\in M_1}\widetilde{A}(\eta\cup\{z\})\longrightarrow\\\prod_{z_1,z_2\in M_0\atop\{z_1,z_2\}\in{X^1}}\widetilde{A}(\eta\cup\{z_1,z_2\})\longrightarrow\prod_{z_1,z_2,z_3\in M_0\atop\{z_1,z_2,z_3\}\in{X^2}}\widetilde{A}(\eta\cup\{z_1,z_2,z_3\})\longrightarrow\ldots\end{array}\right\}$$then we have a short exact sequence$$0\longrightarrow K((y,\widehat{\eta}),M_1')[-1]\longrightarrow K^y(\widehat{\eta},M_0)\longrightarrow K(\widehat{\eta},M_1)\longrightarrow0.$$By induction hypothesis the complexes $K((y,\widehat{\eta}),M_1')$ and $K(\widehat{\eta},M_1)$ are exact; hence the complex $K^y(\widehat{\eta},M_0)$ is exact. The exactness of $K^y(\widehat{\eta},M_0)$ shows the exactness of $K(\widehat{\eta},M_0)$ except at the first non-trivial degree, but at the first non-trivial degree we have already seen exactness. Hence $K(\widehat{\eta},M_0)$ is exact and ${\bf (a)_n}$ is proven.

If we define the complex$$N^y(\widehat{\eta},M_0)=\{0\longrightarrow\frac{\bigcap_{z\in M_1}I(\eta\cup\{z\})}{\bigcap_{z\in M_0}I(\eta\cup\{z\})}\longrightarrow\frac{\bigcap_{z\in M_1'}I(\eta\cup\{y,z\})}{I(\eta\cup\{y\})}\longrightarrow0\longrightarrow\ldots\}$$then we have a short exact sequence$$0\longrightarrow N^y(\widehat{\eta},M_0)\longrightarrow K(\widehat{\eta},M_0)\longrightarrow K^y(\widehat{\eta},M_0)\longrightarrow0.\label{kocut}$$Since we have seen exactness of $K(\widehat{\eta},M_0)$ and of $K^y(\widehat{\eta},M_0)$ we get exactness of $N^y(\widehat{\eta},M_0)$. Now to prove ${\bf (b)_n}$ and ${\bf (c)_n}$ we first suppose $R=\mathbb{Z}$. By exactness of $N^y(\widehat{\eta},M_0)$ we get\begin{gather}\rk_{\mathbb{Z}}(\frac{\widetilde{E}}{ \bigcap_{z\in M_0}I(\eta\cup\{z\})})-\rk_{\mathbb{Z}}(\frac{\widetilde{E}}{ \bigcap_{z\in M_1}I(\eta\cup\{z\})})\notag\\=\rk_{\mathbb{Z}}(\widetilde{A}(\eta\cup\{y\}))-\rk_{\mathbb{Z}}(\frac{\widetilde{E}}{ \bigcap_{z\in M_1'}I(\eta\cup\{y,z\})}).\label{radev}\end{gather}

Observe that the collection $(\prec_{z})_{z\in M_1}$, resp. $(\prec_{z'})_{z'\in M_1'}$, resp. $\prec_y$ is adapted to $(M_1,\widehat{\eta})$, resp. to $(M_1',(y,\widehat{\eta}))$, resp. $(\{y\},\widehat{\eta})$. We associate the sets $G(\widehat{\eta};M_1)$, resp. $G((y,\widehat{\eta});M_1')$, resp. $G(\widehat{\eta};\{y\})$ as before and claim\begin{gather}G(\widehat{\eta};M_0)=G(\widehat{\eta};M_1)\bigcup G(\widehat{\eta};\{y\})\label{cupba}\\G((y,\widehat{\eta});M_1')=G(\widehat{\eta};M_1)\bigcap G(\widehat{\eta};\{y\}).\label{schniba}\end{gather}Here (\ref{cupba}) and $\subset$ in (\ref{schniba}) are very easy. To prove $\supset$ in (\ref{schniba}), let $e_S\in G(\widehat{\eta};M_1)\bigcap G(\widehat{\eta};\{y\})$ and let $e\in S$. Then $e$ is 
$(S,(y,\widehat{\eta}))$-special and $(S,(z,\widehat{\eta}))$-special for some $z\in M_1$. Let $z'\in M_1'$ be the element with $L_y\cap L_z=L_{z'}$. We will show that $e$ is $(S,(z,y,\widehat{\eta}))$-special. If $e\in L_{y}$ then$$e=\max_{\prec_z}({\mathcal{A}}\cap L((z,\widehat{\eta}),\prec_z,S,e))=\max_{\prec_{z'}}({\mathcal{A}}\cap L((z',y,\widehat{\eta}),\prec_{z'},S,e))$$where the first equality follows from the $(S,(z,\widehat{\eta}))$-speciality of $e$ and the second one from $L_y\cap L_z=L_{z'}$. If however $e\notin L_{y}$ then $$e=\max_{\prec_y}({\mathcal{A}}\cap L((y,\widehat{\eta}),\prec_y,S,e))=\max_{\prec_{z'}}({\mathcal{A}}\cap L((z',y,\widehat{\eta}),\prec_{z'},S,e))$$where the first equality follows from the $(S,(y,\widehat{\eta}))$-speciality of $e$ and the second one is clear. 

From (\ref{cupba}) and (\ref{schniba}) we deduce \begin{gather}|G(\widehat{\eta};M_0)|=|G(\widehat{\eta};M_1)|+|G(\widehat{\eta};\{y\})|-|G((y,\widehat{\eta});M_1')|.\label{rang}\end{gather}By induction hypothesis ${\bf (c)_{n-1}}$ we know \begin{align}|G(\widehat{\eta};M_1)| & =\rk_{\mathbb{Z}}(\frac{\widetilde{E}}{ \bigcap_{z\in M_1}I(\eta\cup\{z\})})\notag\\|G(\widehat{\eta};\{y\})| & =\rk_{\mathbb{Z}}(\widetilde{A}(\eta\cup\{y\}))\notag\\|G((y,\widehat{\eta});M_1')| & =\rk_{\mathbb{Z}}(\frac{\widetilde{E}}{ \bigcap_{z\in M_1'}I(\eta\cup\{y,z\})}).\notag\end{align}Thus we may compare (\ref{radev}) with (\ref{rang}) to obtain$$\rk_{\mathbb{Z}}(\frac{\widetilde{E}}{ \bigcap_{z\in M_0}I(\eta\cup\{z\})})=|G(\widehat{\eta};M_0)|.$$Comparing with \ref{erzgrad} we see that source and target of the canonical surjection$$\frac{\widetilde{E}}{J(\widehat{\eta};M_0)}\longrightarrow\frac{\widetilde{E}}{\bigcap_{z\in M_0}I(\eta\cup\{z\})}$$have the same finite $\mathbb{Z}$-rank. But the source is free over $\mathbb{Z}$ (because the $\mathbb{Z}$-submodule $J(\widehat{\eta};M_0)$ of $\widetilde{E}$ has a set of generators each of which is a linear combination, {\it with coefficients in $\{-1,1\}$}, of elements of the obvious (countable) $\mathbb{Z}$-basis of $\widetilde{E}$). Hence this surjection is bijective, so ${\bf (b)_n}$ and ${\bf (c)_n}$ follow if $R=\mathbb{Z}$, and then for arbitrary $R$ by base change.\\

\begin{satz}\label{gilno} For any ${\mathcal A}\subset{\mathbb{P}}(K^{d+1})$, possibly infinite, the ccs $\widetilde{A}$ and $A$ satisfy ${\mathcal{S}}(k)$ for any $1\le k\le d$.
\end{satz}

{\sc Proof:} The condition ${\mathcal{S}}(k)$ for $\widetilde{A}$ requires that for any pointed $(k-1)$-simplex $\widehat{\eta}$ and any subset $M_0\subset N_{\widehat{\eta}}$ which is stable with respect to $\widehat{\eta}$ the sequence$$\widetilde{A}(\eta)\longrightarrow\prod_{z\in {M}_0}\widetilde{A}(\eta\cup\{z\})\longrightarrow\prod_{z_1,z_2\in {M_0}\atop\{z_1,z_2\}\in{X^1}}\widetilde{A}(\eta\cup\{z_1,z_2\})$$is exact. For fixed $\widehat{\eta}$ we may pass to a suitable finite subset of ${\mathcal{A}}$ without changing any of the involved groups. Hence we are in the situation considered above and what we need to show is precisely the exactness of $K(\widehat{\eta},M_0)$ at its first non trivial degree. This we did in \ref{simind}. Hence $\widetilde{A}$ satisfies ${\mathcal{S}}(k)$. But then also $A$ satisfies ${\mathcal{S}}(k)$ because of the split exact sequence (\ref{aspli}).\\

\begin{kor}\label{acyhyp} The ccs $\widetilde{A}$ and $A$ on $X$ are acyclic in positive degrees: for any $k\ge 1$ we have $H^k(X,A)=0$ and $H^k(X,\widetilde{A})=0$.
\end{kor}

{\sc Proof:} \ref{acykrit} and \ref{gilno}.\\ 

For the rest of this section we assume $\kara(K)=0$ and take ${\mathcal{A}}=\mathbb{P}(K^{d+1})$. We write $A_R$ instead of $A$ in order to specify the chosen base ring $R$. Let $\Omega_K^{(d+1)}$ be Drinfel'd's symmetric space of dimension $d$ over $K$. This is the $K$-rigid space obtained by removing all $K$-rational hyperplanes from projective $d$-space ${\mathbb{P}}^{d}_K$. There is a natural ${\rm GL}_{d+1}(K)$-equivariant {\it reduction map} $$r:\Omega_K^{(d+1)}\longrightarrow X$$(see e.g. \cite{ds} for the precise meaning of $r$). For a simplex $\sigma$ of $X$ let $]{\sigma}[=r^{-1}(Star(\sigma))$, the preimage in $\Omega_K^{(d+1)}$ of the star of $\sigma$: the star of $\sigma$ is the union of the open simplices whose closure contains $\sigma$. This $]{\sigma}[$ is an admissible open subset of $\Omega_K^{(d+1)}$, and the collection of all the $]{\sigma}[$ forms an admissible open covering of $\Omega_K^{(d+1)}$. 

\begin{pro}\label{resid} (de Shalit) \cite{ds} For a simplex $\sigma$ of $X$ denote by $H^k_{dR}(]{\sigma}[)$ the $k$-th de Rham cohomology group of the $K$-rigid space $]{\sigma}[$. There is a natural isomorphism$$H^k_{dR}(]{\sigma}[)\cong A_K^k(\sigma).$$
\end{pro}
 
\begin{kor}\label{deraac} (1) (Local acyclicity) Let $\sigma$ be a simplex. For any $k\ge0$ the sequence$$0\longrightarrow H_{dR}^k(\bigcup_{x\in\sigma}]x[)\longrightarrow \prod_{x\in \sigma}H_{dR}^k(]x[)\longrightarrow\prod_{\tau\in X^1\atop\tau\subset\sigma}H_{dR}^k(]\tau[)\longrightarrow\ldots$$is exact.\\(2) (Global acyclicity) (de Shalit) The sequence $$0\longrightarrow H_{dR}^k(\Omega_K^{(d+1)})\longrightarrow\prod_{x\in X^0}H_{dR}^k(]x[)\longrightarrow\prod_{\sigma\in X^1}H_{dR}^k(]\sigma[)\longrightarrow\ldots$$is exact.
\end{kor}

{\sc Proof:} (1) Choose a vertex $x\in\sigma$. Then $M_0=\sigma-\{x\}$ is (as a set of vertices) stable with respect to $x$. Since the ccs $A_K$ satisfies ${\mathcal{S}}(k)$ for any $k$, we derive just as in the proof of \ref{simind} that the sequence$$\prod_{z\in \sigma}A_K(z)\longrightarrow\prod_{\tau\in X^1\atop\tau\subset\sigma}A_K(\tau)\longrightarrow\ldots$$is exact. Inserting \ref{resid} it becomes the exact sequence$$\prod_{x\in \sigma}H_{dR}^k(]x[)\longrightarrow\prod_{\tau\in X^1\atop\tau\subset\sigma}H_{dR}^k(]\tau[)\longrightarrow\ldots.$$On the other hand we have the the spectral sequence$$E_{1}^{rs}=\prod_{\tau\in X^r\atop\tau\subset\sigma}H^s_{dR}(]\tau[)\Longrightarrow H^{s+r}_{dR}(\bigcup_{x\in\sigma}]x[).$$Together (1) follows. The proof of (2) works the same way, using \ref{acyhyp} instead of \ref{gilno}.\\

Corollary \ref{deraac} gives a precise expression of $H_{dR}^k(\Omega_K^{(d+1)})$ through all the $H_{dR}^k(]\sigma[)$. The natural $\mathbb{Z}$-structures $A_{\mathbb{Z}}^k(\sigma)$ in the $A_K^k(\sigma)$ provide natural $\mathbb{Z}$-structures $H_{\mathbb{Z}}^k(]\sigma[)$ in the $H_{dR}^k(]\sigma[)$; hence the ${\rm GL}_{d+1}(K)$-stable subgroup$$H_{\mathbb{Z}}^k(\Omega_K^{(d+1)})=\ke[\prod_{x\in X^0}H_{\mathbb{Z}}^k(]x[)\longrightarrow \prod_{\sigma\in X^1}H_{\mathbb{Z}}^k(]\sigma[)]$$of $H_{dR}^k(\Omega_K^{(d+1)})$. Now our theorem \ref{acyhyp} expresses $H_{\mathbb{Z}}^k(\Omega_K^{(d+1)})$ {\it precisely} through the local terms $A_{\mathbb{Z}}^k(\sigma)$: it tells us that $H_{\mathbb{Z}}^k(\Omega_K^{(d+1)})$ is quasiisomorphic with the complex\begin{gather}\prod_{x\in X^0}H_{\mathbb{Z}}^k(]x[)\longrightarrow \prod_{\sigma\in X^1}H_{\mathbb{Z}}^k(]\sigma[)\longrightarrow\prod_{\sigma\in X^2}H_{\mathbb{Z}}^k(]\sigma[)\longrightarrow\ldots.\label{zrsol}\end{gather}Let us explain why this should have an application to a challenging problem on $p$-adic Abel-Jacobi mappings raised by Raskind and Xarles \cite{raxa}. Let $\Gamma\subset {\rm PGL}_{d+1}(K)$ be a cocompact discrete subgroup such that the quotient $Y=\Gamma\backslash\Omega_K^{(d+1)}$, a smooth projective $K$-scheme, has strictly semistable reduction. For $1\le k\le d$, Raskind and Xarles associate to $Y$ a certain rigid analytic torus $J^k(Y)$, a "$p$-adic intermediate Jacobian"\, ($J^1(Y)$ is the Picard variety of $Y$ and  $J^d(Y)$ is the Albanese variety of $Y$). The device for the construction of $J^k(Y)$ is a canonical $\mathbb{Z}$-structure in the graded pieces $\Gr_*^MH^k(Y)$ of the monodromy filtration on the cohomology $H^k(Y)$ of $Y$ --- both $\ell$-adic ($\ell\ne p$) and log crystalline ($\cong$ de Rham) cohomology. This $\mathbb{Z}$-structure results from the fact that for any component intersection $Z$ of the reduction of $Y$ (so each $Z$ is a smooth projective $k$-scheme), the cycle map$$CH^k(Z\times_k\overline{k})\otimes W(-k)/{\rm tors}\longrightarrow H_{crys}^{2k}(Z\times_k\overline{k}/W)/{\rm tors}$$is bijective (here $W(-k)$ is the ring of Witt vectors $W$ with the action of Frobenius multiplied by $p^k$); similarly for $\ell$-adic cohomology, as was recently proved by Ito \cite{ito}. That is, the $\mathbb{Z}$-structure is essentially given by the collection of Chow groups for all $Z$. Then they define an Abel-Jacobi mapping$$CH^k(Y)_{\rm hom}\longrightarrow J^k(Y)(K)$$with $CH^k(Y)_{\rm hom}$ the group of cycles that are homologically equivalent to zero, using the $\ell$-adic ($\ell\ne p$ and $\ell=p$) Abel-Jabobi mapping which involves the Galois cohomology groups $H_g^1(K,.)$ defined by Bloch and Kato. As they point out, it would be helpful to define the Abel-Jabobi mapping by analytic means.

We expect that such a definition involves $p$-adic integration of cycles on the (contractible !) $K$-rigid space $\Omega_K^{(d+1)}$, similar to Besser's $p$-adic integration on $K$-varieties with good reduction. The link would be the covering spectral sequence $$E_2^{rs}=H^r(\Gamma,H^s(\Omega_K^{(d+1)}))\Longrightarrow H^{r+s}(Y)$$(which exists for both $\ell$-adic ($\ell\ne p$) and de Rham cohomology). Indeed, we know that the associated filtration on $H^{k}(Y)$ {\it is} the monodromy filtration (\cite{desh}, \cite{hk}, \cite{ito}), therefore the ${\mathbb{Z}}$-structure $H_{\mathbb{Z}}^k(\Omega_K^{(d+1)})$ in $H^k(\Omega_K^{(d+1)})$ gives a ${\mathbb{Z}}$-structure in $\Gr_*^MH^k(Y)$. A comparison with that of Raskind and Xarles probably needs the resolution (\ref{zrsol}): the component intersections $Z$ considered by them correspond precisely to the simplices of the quotient simplicial complex $\Gamma\backslash X$.

\section{Local systems arising from representations}

Let $K$, ${\mathcal O}_K$, $\pi$, $k$, $X$ and its orientation be as in section \ref{sechyar}. We fix a natural number $n\ge1$ and let $$U=U^{(n)}=\{g\in {\rm GL}_{d+1}({\mathcal O}_K)|\,\,g\equiv 1\mod\,\,\pi^n\}$$denote the principal congruence subgroup of level $n$ in $G={\rm GL}_{d+1}(K)$. For a vertex $x\in X^0$ we let $$U_x=U_x^{(n)}=gUg^{-1}\,\,\,\mbox{if}\,\,\,x=g([{\mathcal O}_K^{d+1}])\,\,\,\mbox{for some}\,\,\,g\in G$$and for a simplex $\tau=\{x_1,\ldots,x_k\}$ we let$$U_{\tau}=\,\,\mbox{the subgroup generated by}\,\,U_{x_1}\cup\ldots\cup U_{x_k}.$$This is a pro-$p$-group, $p=\kara(k)$.

\begin{lem}\label{ugen} Suppose the lattices $L_z$, $L_{x_1}$ and $L_{x_2}$ represent vertices $z$, $x_1$ and $x_2$ in $X^0$ such that both $x_1$ and $x_2$ are incident to $z$ and such that $L_z=L_{x_1}\cap L_{x_2}$. Then $U_{z}\subset U_{x_1}U_{x_2}$.
\end{lem}

{\sc Proof:} Applying a suitable $g\in G$ we may assume that $L_z={\mathcal O}_K^{d+1}$ and $L_{x_s}=t^s{\mathcal O}_K^{d+1}$ for diagonal matrices $\id\ne t^s=(t_0^s,\ldots,t^s_d)$ satisfying $\{1\}\subset \{t^1_j,t^2_j\}\subset \{1,\pi^{-1}\}$ for all $0\le j\le d$. But then $U_z=U$ as defined above, and $U_{x_s}=t^s U (t^s)^{-1}$ and an easy matrix argument gives the claim.\\

\cite{ssihes} Proposition I.3.1 significantly strengthens \ref{ugen}. It is this interpolation property of the groups $U_x$ which also underlies the acyclicity proof in \cite{ss} and the much more general theory in \cite{ssihes}.

Let $V$ be a smooth representation of $G$ on a (not necessarily free) $\mathbb{Z}[\frac{1}{p}]$-module $V$ which is generated, as a $G$-representation, by its $U$-fixed vectors. Because of $U_{\sigma}\subset U_{\tau}$ if $\sigma\subset \tau$ we can form the hcs $\underline{V}=(V^{U_{\tau}})$ of subspaces of fixed vectors$$V^{U_{\tau}}=\{v\in V|\,\,gv=g\,\,\mbox{for all}\,\,g\in U_{\tau}\}$$with the obvious inclusions as transition maps. In the special case where our $V$ is a $G$ representation on a ${\mathbb C}$-vector space (not just on a $\mathbb{Z}[\frac{1}{p}]$-module), the following theorem (and its version for $n=1$) was proved in \cite{ss}.

\begin{satz}\label{schnstac} Suppose $n>1$. Then the chain complex $C_{\bullet}(X,\underline{V})$ is a resolution of $V$.
\end{satz}

{\sc Proof:} To see $H_k(X,\underline{V})=0$ for $k\ge 1$ it suffices, by \ref{hoacykrit}, to prove ${\mathcal S}^*(k)$, i.e. to prove that for any pointed $(k-1)$-simplex $\widehat{\eta}$ with underlying $(k-1)$-simplex $\eta$ and for any subset $M_0$ of $N_{\widehat{\eta}}$ which is stable with respect to $\widehat{\eta}$, the sequence$$\bigoplus_{z,z'\in {M_0}\atop\{z,z'\}\in{X^1}}V^{U_{\{z,z'\}\cup\eta}}\stackrel{\partial_k}{\longrightarrow} \bigoplus_{z\in {M}_0}V^{U_{\{z\}\cup\eta}}\stackrel{\partial_{k-1}}{\longrightarrow}V^{U_{\eta}}$$is exact. We use induction on $|M_0|$. If $M_0$ is non-empty choose a $y\in M_0$ for which $L_y$ is maximal, i.e. there is no $z\in M_0$ with $L_y\subsetneq L_z$. Then $M_1=M_0-\{y\}$ is stable with respect to $\widehat{\eta}$. Letting $$M_1'=\{z'\in M_1|\,\,\{y,z'\}\in X^1\}$$we first claim that\begin{gather}\bigoplus_{{z'\in {M_1'}}}V^{U_{\{y,z'\}\cup\eta}}\longrightarrow V^{U_y}\bigcap\sum_{z\in M_1}V^{U_{\{z\}\cup\eta}}\label{suss}\end{gather}is surjective. Let $v=\sum_{z\in M_1}v_z$ be an element of the right hand side with $v_z\in V^{U_{\{z\}\cup\eta}}$ for all $z\in M_1$. Since $V$ is smooth we can find a finitely generated submodule $V'$ of $V$ containing all the $v_z$ which is stable under $U_{y}$. The action of $U_y$ on $V'$ factors through a finite quotient of $U_y$; let $\overline{U}_{y}\subset U_y$ be a set of representatives for this quotient. Since $v$ is fixed by ${U}_{y}$ it follows that$$v=\frac{1}{|\overline{U}_{y}|}\sum_{g\in \overline{U}_{y}}g.v=\sum_{z\in M_1}\frac{1}{|\overline{U}_{y}|}\sum_{g\in \overline{U}_{y}}g.v_z.$$Since $M_0$ is stable, there exists for any $z\in M_1$ a $z'\in M_1'$ such that $L_{z'}=L_z\cap L_y$. It will be enough to show $\sum_{g\in\overline{U}_{y}}g.v_z\in V^{U_{\{y,z'\}\cup\eta}}$. The stability under $U_y$ is clear.
 Now let $h\in U_{z'}$. By \ref{ugen} we may factor $h$ as $h=h_yh_z$ with $h_y\in U_y$ and $h_z\in U_{z}$. Since $n>1$ and since there is a vertex incident to both $z$ and $y$ we have $g^{-1}U_{z}g=U_{z}$ for any $g\in{U}_{y}$, hence $h_zg=gh^g_z$ with $h^g_z\in U_{z}$. Thus$$h\sum_{g\in \overline{U}_{y}}g.v_z=h_y\sum_{g\in \overline{U}_{y}}gh^g_z.v_z=\sum_{g\in \overline{U}_{y}}g.v_z,$$i.e. $\sum_{g\in \overline{U}_{y}}g.v_z$ is stable under $U_{z'}$. Finally let $h\in U_x$ for some $x\in\eta$. Since $x$ and $y$ are incident we have $g^{-1}U_{x}g=U_{x}$, hence there is for any $g\in{U}_{y}$ a $h^g\in U_{x}$ with $hg=gh^g$. Then$$h\sum_{g\in \overline{U}_{y}}g.v_z=\sum_{g\in \overline{U}_{y}}gh^g.v_z=\sum_{g\in \overline{U}_{y}}g.v_z$$so we have shown stability under $U_{\eta}$. The surjectivity of (\ref{suss}) is proven. Now let $c=(c_z)_{z\in M_0}$ be an element of $\ke({\partial_{k-1}})$. Then necessarily$$c_y\in V^{U_y}\bigcap\sum_{z\in M_1}V^{U_{\{z\}\cup\eta}}.$$By the surjectivity of (\ref{suss}) we may therefore modify $c$ by an element of $\bi(\partial_k)$ such that for the new $c=(c_z)_{z\in M_0}\in\ke(\partial_{k-1})$ we have $c_y=0$. But then the induction hypothesis, applied to $M_1$, tells us that after another such modification we can achieve $c=0$. We have shown that $C_{\bullet}(X,\underline{V})$ is exact in positive degrees. It remains to observe that the hypothesis that $V$ is generated by $V^U$ is equivalent with the surjectivity of$$C_{0}(X,\underline{V})=\bigoplus_{x\in X^0}V^{U_x}\longrightarrow V.$$



\begin{thebibliography}{abcdefgh}  
\bibitem{alon}{\it G. Alon}, Cohomology of local systems coming from $p$-adic hyperplane arrangements, to appear in Israel J. Math.
\bibitem{bro} {\it K. S. Brown}, Buildings, Berlin-Heidelberg-New York: Springer (1989) 
\bibitem{ds} {\it E. de Shalit}, Residues on buildings and de Rham cohomology of $p$-adic symmetric domains, Duke Math. J. {\bf 106} (2001), no.1, 123--191
\bibitem{desh}{\it E. de Shalit}, The $p$-adic monodromy-weight conjecture for $p$-adically uniformized varieties, Compositio Math. {\bf 141} (2005), no. 1, 101--120
\bibitem{hk}{\it E. Grosse-Kl\"onne}, Frobenius and Monodromy operators in rigid analysis, and Drinfel'd's symmetric space, Journal of Algebraic Geometry 14 (2005), 391--437
\bibitem{int}{\it E. Grosse-Kl\"onne}, Integral structures in the $p$-adic holomorphic discrete series, Representation Theory 9 (2005), 354--384
\bibitem{ito}{\it T. Ito}, Weight-Monodromy conjecture for $p$-adically uniformized varieties, preprint 2003
\bibitem{raxa}{\it W. Raskind and X. Xarles}, On the \'{e}tale cohomology of algebraic varieties with totally degenerate reduction over $p$-adic fields, preprint 2003 
\bibitem{schn}{\it P. Schneider}, The cohomology of local systems on 
$p$-adically uniformized varieties, Math. Ann. {\bf 293}, 623--650 (1992)
\bibitem{ss}{\it P. Schneider and U. Stuhler}, Resolutions for smooth representations of the general linear group over a local field, J. reine angew. Math. {\bf 436}, 19--32 (1993)
\bibitem{ssihes}{\it P. Schneider and U. Stuhler}, Representation theory and sheaves on the Bruhat-Tits building, Publ. Math. IHES {\bf 85}, 97 -- 191 (1997)
\end{thebibliography}
\end{document}